\documentclass[12pt,twoside]{amsart}
\usepackage{amssymb}
\usepackage{amscd}

\title{Applications of Kawamata's positivity theorem}  
\author{Osamu Fujino} 
\date{\today}
\address{Research Institute for Mathematical Sciences\\ 
 Kyoto University, Kyoto 606-8502 Japan}
\email{fujino@kurims.kyoto-u.ac.jp}

\newcommand{\bQ}[0]{{\mathbb Q}}

\newcommand{\Supp}[0]{{\operatorname{Supp}}}
\newcommand{\codim}[0]{{\operatorname{codim}}}
\newcommand{\NE}[0]{{\operatorname{NE}}}

\newtheorem{thm}{Theorem}[section]

\newtheorem{cor}[thm]{Corollary}

\theoremstyle{definition}
\newtheorem{defn}[thm]{Definition}
\newtheorem{assume}[thm]{Assumptions}
\newtheorem{prob}[thm]{Problem}
\newtheorem{rem}[thm]{Remark}
\newtheorem{ack}{Acknowledgments}        

\theoremstyle{remark}
\newtheorem*{claim}{Claim}

\begin{document}
\bibliographystyle{amsalpha+}
\maketitle

\abstract In this paper we treat some applications of Kawamata's positivity 
theorem. 
We get a weak answer to \cite [Section 3]{KeMaMc}. 
And we investigate the singularities on the target spaces of some 
morphisms. 
\endabstract

\setcounter{section}{-1}
\section{Introduction}\label{se1}

In this paper we treat some applications of Kawamata's positivity 
theorem (See Theorem (\ref {kpt})), 
which are related to the following important problem. 

\begin{prob}\label{qu}
Let $(X,\Delta)$ be a proper klt pair. 
Let $f:X\to S$ be a proper surjective morphism 
onto a normal variety $S$ with connected fibers. 
Assume that $K_X+\Delta\sim_{\bQ,f}0$. 
Then is there any effective $\bQ$-divisor $B$ on $S$ 
such that 
$$
K_X+\Delta\sim _{\bQ} f^{*} (K_S+B)
$$ 
and that the pair $(S,B)$ is again klt? 
\end{prob}

A special case of 
Problem (\ref {qu}) was studied in \cite [Section 3]{KeMaMc}, 
where general fibers are rational curves 
as a step toward the proof of the three dimensional log abundance theorem. 
Thanks to Kodaira's canonical bundle formula 
for elliptic surfaces and \cite [Section 3]{KeMaMc}, 
the problem is affirmative under the assumption that $\dim X \leq 2$. 
However, the problem is much harder in higher dimension 
because of the lack of canonical bundle formula and so forth. 
In this paper, we prove the following theorem as an application 
of positivity theorem of Kawamata, which could be 
viewed as a partial answer to Problem (\ref{qu}). 
This is the main theorem of this paper. 

\begin{thm}\label{mt}
Let $(X,\Delta)$ be a proper sub klt pair. 
Let $f:X\to S$ be a proper surjective morphism 
onto a normal variety $S$ with connected fibers. 
Assume that 
$\dim _{k(\eta)} f_{*}{\mathcal O}_X (\ulcorner -\Delta\urcorner)
\otimes_{{\mathcal O}_{S}} k(\eta)=1$, 
where $\eta$ is the generic point of $S$. 
And assume that $K_X+\Delta\sim_{\bQ,f}0$, that is, there exists a 
$\bQ$-Cartier $\bQ$-divisor $A$ 
on $S$ such that $K_X+\Delta\sim_{\bQ}f^{*} A$. 
Let $H$ be an ample Cartier divisor on $S$, and $\epsilon$ a 
positive rational number. 
Then there exists a $\bQ$-divisor $B$ on $S$ 
such that 
$$
K_S+B\sim_{\bQ}A+\epsilon H ,
$$
$$
K_X+\Delta+\epsilon f^{*} H\sim _{\bQ} f^{*} (K_S+B),
$$ 
and that the pair $(S,B)$ is sub klt. 

Furthermore, if 
$f_{*}{\mathcal O} _{X} (\ulcorner -\Delta\urcorner) ={\mathcal O}
_S$, 
then we can make $B$ effective, that is, $(S,B)$ is klt. 
In particular, $S$ has only rational singularities. 
\end{thm}

By using Theorem (\ref {mt}), 
we obtain the cone theorem for the base space $S$ (See Theorem (\ref {CCC})), 
whose proof is given in Section \ref {c} for the reader's convenience. 
We also prove that the target space of an extremal contraction 
is at worst ``Kawamata log terminal'' 
(See Corollary (\ref {cor})). 
Corollary (\ref{nakayama}) is a reformulation of 
a result of Nakayama. 

In this paper, we will work over $\mathbb C$, the complex number 
field, and make use of the standard notations as in \cite [Notation 0.4] 
{KM}. 

\begin{ack}
I would like to express my sincere 
gratitude to Professors Shigefumi Mori and 
Noboru Nakayama for many useful comments. 
I am grateful to Dr.~Daisuke Matsushita, 
who informed me of the paper 
\cite {N1}, which was a starting 
point of this paper. 
I am also grateful to Professor Yoichi Miyaoka 
for warm encouragements.   
\end{ack}

\section{Definitions and Preliminaries}

We make some definitions and cite the key theorem in this section.

\begin{defn}\label{def1}
Let $f:X\to S $ be a proper surjective morphism of normal varieties 
with connected fibers. 
\begin{enumerate}
\item [(i)] A divisor $D$ is called 
$f$-exceptional if $\codim _S f(D)\geq 2$. 
\item [(ii)] Two $\bQ$-divisors $\Delta$ and $\Delta'$ on $X$ are called 
$\bQ$-linearly $f$-equivalent, denoted by $\Delta\sim _{\bQ,f} \Delta'$, 
if there exists 
a positive integer $r$ such that $r\Delta$ and $r\Delta'$ are linearly 
$f$-equivalent (See \cite [Notation 0.4 (5)]{KM}).
\end{enumerate}
\end{defn}

\begin{defn}
A pair $(X,\Delta)$ of normal variety and a 
$\bQ$-divisor $\Delta=\sum_i d_i \Delta_i$ 
is said to be sub Kawamata log terminal (sub klt, for short) 
(resp. ~divisorial log terminal (dlt, for short)) if 
the following conditions are satisfied: 
\begin{enumerate}
\item [(1)] $K_X+\Delta$ is a $\bQ$-Cartier $\bQ$-divisor ;
\item [(2)] $d_i < 1$ (resp. ~$0\leq d_i \leq 1$) ;
\item [(3)] there exists a log resolution (See \cite [Notation 0.4 (10)] 
{KM}) $\mu : Y \to X$ such that $a_j>-1$ for all $j$ 
in the canonical bundle formula,  
$$
K_Y+\mu_{*}^{-1} \Delta=\mu^{*}(K_X+\Delta)+\sum _j a_j E_j.
$$
\end{enumerate}
A pair $(X,\Delta)$ is called Kawamata log terminal (klt, for short) 
if $(X,\Delta)$ is sub klt and $\Delta$ is effective. 
We say that a variety $X$ has only canonical 
singularities if 
$(X,0)$ is klt and $a_j\geq 0$ for all $j$ in (3). 
\end{defn}

The notion of divisorial log terminal pair was first 
introduced by V.~V.~Shokurov in his paper \cite {Sh} 
(for another equivalent definition, see \cite 
[Divisorial Log Terminal Theorem]{S} and \cite [Definition 2.37]{KM}). 

\begin{defn}
Let $f:X\to S $ be a smooth surjective morphism of varieties 
with connected fibers. 
A reduced effective divisor $D=\sum _i D_i$ 
on $X$ such that $D_i$ is mapped onto $S$ 
for every $i$ is said to be relatively normal crossing 
if the following condition holds. 
For each closed point $x$ of $X$, 
there exists an open neighborhood $U$ 
(with respect to the classical topology) 
and 
$u_1, \cdots, u_k \in {\mathcal O}_{X,x}$ 
inducing a regular system of parameters 
on $f^{-1} f(x)$ at $x$, 
where $k=\dim _x f^{-1}f(x)$, 
such that 
$D\cap U=\{ u_1\cdots u_l=0\}$ for some $l$ such that $0\leq l\leq k$. 
\end{defn}

The next theorem is \cite [Theorem 2]{Ka3}, which plays an essential role in 
this paper. 
The conditions (2) and (3) are different from the original ones. 
But we do not have to change the proof in \cite {Ka3}. 

\begin{thm}[Kawamata's positivity theorem]\label{kpt}
Let $g: Y\to T$ be a surjective morphism of smooth projective varieties with
connected fibers.
Let $P = \sum _{j} P_j$ and $Q = \sum_{l} Q_{l}$ be
normal crossing divisors on $Y$ and $T$, respectively, such that
$g^{-1}(Q) \subset P$ and 
$g$ is smooth over $T \setminus Q$.
Let $D = \sum_j d_jP_j$ be a $\bQ$-divisor on $Y$ 
{\em{(}}$d_j$'s may be negative{\em{)}}, 
which satisfies the following conditions:
\begin{enumerate}
\item [(1)] $D = D^h + D^v$ such that 
every irreducible component of $D^h$ is mapped surjecively onto $T$ by $g$, 
$g: \Supp(D^h) \to T$ is relatively normal crossing 
over $T \setminus Q$, and
$g(\Supp(D^v)) \subset Q$.  An irreducible 
component of $D^h$ (resp. ~$D^v$) is called 
{\it horizontal} (resp. {\it vertical}).
\item [(2)] $d_j < 1$ if $P_j$ is not $g$-exceptional.
\item [(3)]
$\dim _{k(\eta)} g_{*}{\mathcal O}_Y (\ulcorner -D\urcorner)
\otimes_{{\mathcal O}_{T}} k(\eta)=1$, 
where $\eta$ is the generic point of $T$. 
\item [(4)] $K_Y + D \sim_{\bQ} g^*(K_T + L)$ for some 
$\bQ$-divisor $L$ on $T$.
\end{enumerate}
Let 
\begin{eqnarray*}
g^{*}Q_{l} &=&\sum_j w_{l j}P_j \\
\bar d_j &=& \frac {d_j + w_{l j} - 1}{w_{l j}} 
\text{ if } g(P_j) = Q_{l} \\  
\delta_{l} &=& \text{max }\{ \bar d_j; g(P_j) = Q_{l}\} \\
\Delta _{0}&=& \sum_{l} \delta_{l}Q_{l} \\
M &= &L - \Delta_{0}.
\end{eqnarray*}
Then $M$ is nef.
\end{thm}

\section{Proof of the Main Theorem}

\proof [Proof of Theorem (\ref {mt})](cf. \cite [Theorem 2]{N1}) 
By using the desingularization theorem (cf. \cite [Resolution Lemma]{S}) we 
have the following commutative diagram:
$$
\begin{CD}
Y @>\text{$\nu$}>> X \\
  @V\text {$g$}VV  @VV\text {$f$}V\ \ , \\
T @>\text{$\mu$}>>S
\end{CD}
$$
where 
\begin{enumerate}
\item [(i)] $Y$ and $T$ are smooth projective varieties, 
\item [(ii)] $\nu$ and $\mu$ are projective birational morphisms, 
\item [(iii)] we define $\bQ$-divisors $D$ and $L$ on $Y$ and $T$ 
by the following relations:
$$
K_Y+D=\nu^{*}(K_X+\Delta),
$$
$$
K_T+L  \sim _{\bQ} \mu^{*} A,
$$
\item [(iv)] there are simple normal crossing divisors $P$ and $Q$ on 
$Y$ and $T$ such that they satisfy the conditions of Theorem (\ref {kpt}) and 
there exists a set of positive rational numbers 
$\{s_l\}$ such that $\mu^{*} H-\sum_{l} s_l Q_l$ is ample. 
\end{enumerate}
By the construction, the conditions (1) and (4) of Theorem 
(\ref {kpt}) are satisfied. 
Since $(X,\Delta)$ is sub klt, the condition (2) of Theorem 
(\ref {kpt}) is satisfied. 
The condition (3) of Theorem 
(\ref {kpt}) can be checked by the following claim. 
Note that $\mu$ is birational. 
We put $h:=f\circ \nu$.
\begin{claim}[A]\label{A}
${\mathcal O}_{S} \subset h_{*}{\mathcal O}_{Y}(\ulcorner -D\urcorner) 
\subset f_{*}{\mathcal O}_{X}(\ulcorner -\Delta\urcorner)$.
\end{claim}

\proof[Proof of Claim (A)]
Since ${\mathcal O}_{Y}\subset {\mathcal O}_{Y}(\ulcorner -D\urcorner)$, 
we have  
${\mathcal O} _{S} 
=h_{*}{\mathcal O} _{Y}
\subset h_{*}{\mathcal O}_{Y}(\ulcorner -D\urcorner)$.
Note that $\ulcorner -D\urcorner=\nu _{*}^{-1}
\ulcorner -\Delta\urcorner +F$, where $F$ is effective and $\nu$-exceptional. 
Then 
$$
\begin{array}{clcl}
&\Gamma(U, \nu_{*}{\mathcal O}_{Y}(\ulcorner -D\urcorner))&&\\
\subset&  \Gamma(U \backslash 
\nu(F), \nu_{*}{\mathcal O}_{Y}(\ulcorner -D\urcorner))
&=&\Gamma(U \backslash 
\nu(F), \nu_{*}{\mathcal O}_{Y}(\nu_{*}^{-1}\ulcorner -\Delta\urcorner))\\
\subset&  \Gamma(U \backslash 
\nu(F), {\mathcal O}_{X}(\ulcorner -\Delta\urcorner))
&=&\Gamma(U, {\mathcal O}_{X}(\ulcorner -\Delta\urcorner)),
\end{array}
$$
where $U$ is a Zariski open set of $X$. 
So we have $\nu_{*}{\mathcal O}_{Y}(\ulcorner -D\urcorner)\subset {\mathcal O}
_{X}(\ulcorner -\Delta\urcorner)$. 
Then $h_{*}{\mathcal O}_{Y}(\ulcorner -D\urcorner) 
\subset f_{*}{\mathcal O}_{X}(\ulcorner -\Delta\urcorner)$. 
We get Claim (A). 
\endproof

So we can apply Theorem (\ref {kpt}) to $g:Y\to T$. 
The divisors $\Delta_{0}$ and $M$ are as in Theorem (\ref {kpt}). 
Then $M$ is nef. 
Since $M$ is nef, we have that 
$$
M+\epsilon\mu^{*}H-\epsilon'\sum_{l} s_l Q_l 
$$ 
is ample for $0<\epsilon'\leq \epsilon$. 
We take a general Cartier divisor 
$$ 
F_0\in |m(M+\epsilon\mu^{*}H-\epsilon'\sum_{l} s_l Q_l )|
$$
for a sufficiently large and divisible integer $m$. 
We can assume that $\Supp (F_0 \cup \sum _l Q_l)$ is a simple normal 
crossing divisor. 
And we define $F:=(1/m)F_{0}$. 
Then 
$$
L+\epsilon\mu^{*}H\sim_{\bQ} F+\Delta_0 + \epsilon'\sum_{l} s_l Q_l .
$$
Let $B_0:= F+\Delta_0 + \epsilon'\sum_{l} s_l Q_l$ and 
$\mu_* B_0=B$. 
We have $K_T+B_{0}=\mu ^{*} (K_S+B)$. 
By the definition, $\llcorner\Delta_0\lrcorner\leq 0 $. 
So $\llcorner F+\Delta_0 + \epsilon'\sum_{l} s_l Q_l \lrcorner \leq 0$ 
when $\epsilon'$ is small enough. 
Then $(S,B)$ is sub klt. 
By the construction we have  
$$
K_S+B\sim_{\bQ}A+\epsilon H ,
$$
$$
K_X+\Delta+\epsilon f^{*} H\sim_{\bQ} f^{*}(K_S+B).
$$
If we assume furthermore that 
$f_{*}{\mathcal O}_X (\ulcorner -\Delta\urcorner)={\mathcal O}_S$, 
we can prove the following claim. 

\begin{claim}[B]
If $\mu_* Q_l\ne 0$, then $\delta_l\geq 0$.
\end{claim}
\proof[Proof of Claim (B)]
If $\ulcorner -d_{j}\urcorner\geq w_{lj}$ for every $j$, 
then $\ulcorner -D\urcorner\geq g^{*}Q_{l}$. 
So $g_{*} {\mathcal O} _{Y}( \ulcorner -D\urcorner) \supset 
{\mathcal O} _{T} (Q_l)$. 
Then $\mathcal O _{S} =h_{*}{\mathcal O} _{Y}( \ulcorner -D\urcorner)
\supset \mu _{*} {\mathcal O} _{T} (Q_l)$ by Claim (A). 
It is a contradiction. 
So we have that $\ulcorner -d_j\urcorner <w_{lj}$ for some $j$. 
Since $w_{lj}$ is an integer, we have that 
$-d_j+1\leq w_{lj}$. 
Then $\bar d_j\geq 0$. We get $\delta_l\geq 0$. 
\endproof

So $B$ is effective if  
$f_{*}{\mathcal O}_X (\ulcorner -\Delta\urcorner)={\mathcal O}_S$. 
This completes the proof. 
\endproof

Note that Theorem (\ref {mt}) implies a generalization of Koll\'ar's 
result \cite [Remark 3.16]{Ko2}. 

\section{Applications of the Main Theorem}

The following theorem is the cone theorem for $(S,A-K_S)$. 
This implies the argument in \cite [(5.4.2)]{C}.

\begin{thm}[Generalized Cone Theorem]\label{CCC}
In Theorem (\ref{mt}) we assume that 
$f_{*}{\mathcal O}_X (\ulcorner -\Delta\urcorner)={\mathcal O}_S$. 
Then we have the cone 
theorem of $S$ as follows: 
\begin{enumerate}
\item [(1)] There are {\em{(}}possibly 
countably  many{\em{)}} rational curves $C_j\subset S$ 
such that $\mathbb R_{\geq 0}[C_j]$ is an $A$-negative extremal 
ray for every $j$ and 
$$
\overline{\NE}(S)=
\overline{\NE}(S)_{A\geq 0}+\sum_j \mathbb R _{\geq 0} [C_j].
$$
\item [(2)] For any $\delta>0$ and every ample $\bQ$-divisor $F$, 
$$
\overline{\NE}(S)=\overline{\NE}(S)_{(A+\delta F)\geq 0}+
\sum_{\text{\rm {finite}}} \mathbb R _{\geq 0} [C_j].
$$
\item [(3)] The contraction theorem holds for any $A$-negative extremal face 
{\em(}for more precise statement, 
see \cite [Theorem 3.7 (3), (4)]{KM}{\em)}.
\end{enumerate}
\end{thm}

\proof
See Section \ref {c}.
\endproof

The next theorem is a partial answer to Problem (\ref {qu}) 
under some assumptions.  

\begin{thm}\label{aa}
Let $(X,\Delta)$ be a proper sub klt pair. 
Let $f:X\to S$ be a proper surjective morphism 
onto a normal projective variety with connected fibers. 
Assume that $K_X+\Delta\sim_{\bQ,f}0$ 
and $f_{*}{\mathcal O} _{X} (\ulcorner -\Delta\urcorner) ={\mathcal O}_S$. 
Assume that 
$S$ is $\bQ$-factorial and the Picard number $\rho (S)=1$, 
and the irregularity $q(S)=0$. 
Then there is an effective $\bQ$-divisor $\Delta'$ on $S$ 
such that 
\begin{gather}
\tag{!} K_X+\Delta\sim _{\bQ} f^{*} (K_S+\Delta'),
\end{gather}
and that the pair $(S,\Delta')$ is klt. 
\end{thm}

\proof
We use the same notations as in the proof of Theorem (\ref {mt}). 
By Theorem (\ref{mt}) and the $\bQ$-factoriality of $S$, 
we have that $(S,\mu_{*}\Delta_0)$ is klt. 
The $\bQ$-divisor $\mu_{*}M$ is an ample $\bQ$-Cartier $\bQ$-divisor 
or $\mu_{*}M\sim_{\bQ}0$ since $S$ is 
$\bQ$-factorial and $\rho (S)=1$, and $q(S)=0$. 
When $\mu_{*}M\sim_{\bQ}0$, we put $\Delta':=\mu_{*}\Delta_0$. 
So $(S,\Delta')$ is klt and satisfies (!). 
When $\mu_{*}M$ is ample, we take a sufficiently large and divisible 
integer $k$ such that $|k\mu_{*}M|$ is very ample. 
Let $C$ be a general member of $|k\mu_{*}M|$. 
We put $\Delta':=(1/k)C+\mu_{*}\Delta_0$. 
Then $(S,\Delta')$ is klt and satisfies (!). 
\endproof

\begin{rem}
In Theorem (\ref {aa}), the assumption $q(S)=0$ is satisfied if 
$-K_S$ is nef and big. 
It is because $S$ is klt by Theorem (\ref {mt}) and the 
$\bQ$-factoriality of $S$. 
So $q(S)=h^{1}(S,{\mathcal O}_S)=0$ by Kawamata-Viehweg 
Vanishing Theorem.
\end{rem}

\begin{rem}
On the assumption that $X$ has only canonical singularities 
and the general fibers of $f$ are smooth elliptic curves, 
Problem (\ref {qu}) was proved (See \cite [Corollary 0.4]{N0}). 
\end{rem}

The next corollary is the generalization of Koll\'ar's theorem 
(See \cite [Corollary 7.4]{Ko1}). 
\begin{cor}\label{cor}
Let $(X,\Delta)$ be a projective dlt pair. 
Let $f:X\to S$ be an extremal contraction 
{\em{(}}See \cite [Theorem 3-2-1]{KMM}{\em{)}}. 
Then there exists an effective $\bQ$-divisor $\Delta'$ on $S$ 
such that $(S,\Delta')$ is klt. 
In particular, $S$ has only rational singularities. 
\end{cor}

\proof
Let $H$ be an ample $\bQ$-Cartier $\bQ$-divisor on $S$ 
such that $H':=-(K_X+\Delta)+f^{*}H$ is ample. 
Let $m$ be a positive integer such that $m\Delta$ is a 
$\mathbb Z$-divisor. 
Let $m'$ be a sufficiently large and divisible integer such that 
$\mathcal O_X(m\Delta+m'H')$ is generated by global sections. 
We take a general member $D'\in |m\Delta+m'H'|$. 
Then $H'\sim _{\bQ} (1/{m'})(D'-m\Delta)$. 
So we have that 
$K_X+\Delta+(\epsilon/{m'}) (D'-m\Delta)$ is $\bQ$-Cartier and klt 
for any rational number $0<\epsilon\ll 1$ (See \cite [Proposition 2.43] {KM}). 
We take a sufficiently large and divisible integer $k$ such that $kH'$ 
is very ample. 
Let $D''$ be a general member of $|kH'|$. 
Then 
$$ 
(X,\Delta+\frac{\epsilon}{m'}(D'-m\Delta) +\frac{1-\epsilon}{k} D'') 
$$
is klt and 
$$
K_X+\Delta+\frac{\epsilon}{m'}(D'-m\Delta) +\frac{1-\epsilon}{k} D''
\sim _{\bQ,f} 0.
$$
Apply Theorem (\ref {mt}) for 
$$
f: (X, 
\Delta+\frac{\epsilon}{m'}(D'-m\Delta) +\frac{1-\epsilon}{k} D'')
\to S.
$$
We get the result. 
\endproof

\begin{cor}
Let $f:X\to S$ be a Mori fiber space 
{\em {(}}for the definition of a Mori 
fiber space, see \cite [(1.2)]{C}{\em {)}}. 
Then $S$ is klt. 
\end{cor}

\proof
Apply Corollary (\ref{cor}). Note that $S$ is $\bQ$-factorial 
(See \cite [Corollary 3.18]{KM}). 
\endproof

The following corollary is a reformulation of 
\cite [Corollary A.4.4]{N2}, whose assumption 
is slightly different from ours. 
Our proof is much simpler than \cite [Corollary A.4.4]{N2}, 
but we can only treat the global situation. 
For the non-projective case, 
we refer the reader to \cite [Appendix]{N2}. 

\begin{cor}\label{nakayama}
Let $(X,\Delta)$ be a proper sub klt pair. 
Let $f:X\to S$ be a proper surjective morphism 
onto a normal projective variety $S$ with connected fibers. 
Assume that 
$f_{*}{\mathcal O} _{X} (\ulcorner -\Delta\urcorner) ={\mathcal O}_S$ 
and $-(K_X+\Delta)$ is $f$-nef and $f$-abundant 
{\em{(}}See \cite [Definition 6-1-1]{KMM}{\em{)}}. 
Then there exists an 
effective  $\bQ$-divisor $\Delta'$ on $S$ 
such that $(S,\Delta')$ is klt. 
In particular, $S$ has only rational singularities. 
\end{cor}
\proof
By \cite [Proposition 6-1-3]{KMM}, there exists a diagram: 
$$
\begin{CD}
Y @>\text{$\mu$}>> Z \\
  @V\text {$g$}VV  @VV\text {$\nu$}V \\
X @>\text{$f$}>>S
\end{CD}
$$
which satisfies the following conditions; 
\begin{enumerate}
\item [(i)] it is commutative, that is, $h:=f\circ g= \nu\circ\mu$, 
\item [(ii)] $\mu$, $\nu$ and $g$ are projective morphisms, 
\item [(iii)] $Y$ and $Z$ are nonsingular varieties, 
\item [(iv)] $g$ is a birational morphism and $\mu$ is a surjective morphism 
with connected fibers, and 
\item [(v)] there exists a $\nu$-nef and $\nu$-big $\bQ$-Cartier 
$\bQ$-divisor $D$ on $Z$ such that 
$$
K_{Y}+\Delta'':=g^{*}(K_{X}+\Delta)\sim _{\bQ}\mu^{*} (-D).
$$  
\end{enumerate}
By \cite [Lemma 1.7]{Ka0}, there exists an effective $\bQ$-divisor $C_0$ 
on $Z$ such that $D-C_0$ is $\nu$-ample. 
We define $C:=\mu ^{*} C_0$. 
If $m$ is a sufficiently large integer, then 
$(Y,\Delta''+(1/m)C)$ is sub klt. 
We put 
$$
H:=-(K_Y+\Delta'')-\frac{1}{m} C\sim _{\bQ} \mu ^{*}( D-\frac {1}{m} C_0).
$$
Then $H$ is $h$-semi-ample since 
$D-(1/m)C_0$ is $\nu$-ample. 
We take a very ample Cartier divisor $A$ on $S$ such that 
$H+h^{*} A$ is semi-ample. 
Let $E$ be a general member of 
$|k(H+h^{*}A)|$ for a sufficiently large and divisible integer $k$. 
Then 
$$
(Y,\Delta''+\frac {1}{m} C+ \frac{1}{k} E) 
$$
is sub klt and 
$$
K_Y+\Delta''+\frac {1}{m} C+ \frac{1}{k} E\sim _{\bQ,h}0
$$
by the construction. 
So we can apply Theorem (\ref {mt}) for 
$$
h: (Y,\Delta''+\frac {1}{m} C+ \frac{1}{k} E)\to S. 
$$
Note that 
$$
{\mathcal O}_{S}\subset 
h_{*}{\mathcal O} _{Y}(\ulcorner -\Delta'' -\frac{1}{m}C-\frac{1}{k}E\urcorner)
=h_{*}{\mathcal O} _{Y}(\ulcorner -\Delta''\urcorner)
\subset f_{*}{\mathcal O} _{X}(\ulcorner -\Delta\urcorner)
$$
by Claim (A) in the proof of Theorem (\ref {mt}). 
This complete the proof.
\endproof

\section{Generalized Cone Theorem}\label{c}

In this section we always work on the following assumption. 

\begin{assume}
Let $S$ be a normal projective variety and $A$ a $\bQ$-Cartier $\bQ$-divisor 
on $S$. 
For any positive rational number $\epsilon$ 
and every ample Cartier divisor $H$, 
there exists an effective $\bQ$-divisor $B$ on $S$ such that 
$K_S+B\sim_{\bQ}A+\epsilon H$ and that $(S,B)$ is klt. 
\end{assume}

\begin{defn}
Let $F$ be an ample Cartier divisor. 
We define 
$$
r:=\sup\{t\in \mathbb R ; F+tA \ \text {is nef}\ \}.
$$
\end{defn}

\begin{thm}[Generalized Cone Theorem]\label{cone}
We have the generalization of the cone theorem as follows: 
\begin{enumerate}
\item [(1)] There are {\em{(}}possibly 
countably  many{\em{)}} rational curves $C_j\subset S$ 
such that $\mathbb R_{\geq 0}[C_j]$ is an $A$-negative extremal 
ray for every $j$ and 
$$
\overline{\NE}(S)=
\overline{\NE}(S)_{A\geq 0}+\sum_j \mathbb R _{\geq 0} [C_j].
$$
\item [(2)] For any $\delta>0$ and every ample $\bQ$-divisor $F$, 
$$
\overline{\NE}(S)=\overline{\NE}(S)_{(A+\delta F)\geq 0}+
\sum_{\text{\rm {finite}}} \mathbb R _{\geq 0} [C_j].
$$
\item [(3)] The contraction theorem holds for any $A$-negative extremal face 
{\em(}for more precise statement, 
see \cite [Theorem 3.7 (3), (4)]{KM}{\em)}. 
\end{enumerate}
\end{thm}

\proof 
If $A$ is nef, then there is nothing to be proved. 
So we can assume that $A$ is not nef. 
Then (2) is obvious. 
Note that for rational numbers 
$0<\delta'<\delta$ 
there is an effective 
$\bQ$-divisor $B'$ such that 
$A+\delta F\sim_{\bQ} K_S+B'+\delta'F$ 
and $(S,B')$ is klt. 
So we can reduce it to the well-known cone theorem for klt pairs 
(See \cite[Theorem 3.7]{KM}). 

Let $\epsilon$ be a small positive 
rational number and $K_S+B\sim_{\bQ}A+\epsilon F$. 
If $F+r_0(K_S+B)$ is nef but not ample, then $r_0$ is a rational number 
by the rationality theorem for klt pairs. 
So we get the rationality of $r$. 
By \cite [Lemma 4-2-2]{KMM} and the rationality of $r$, we have 
$$
\overline{\NE}(S)=\overline{\overline{\NE}(S)_{A\geq 0}+
\sum_j \mathbb R _{\geq 0} [C_j]},
$$
where the right hand side is the closure of the cone 
generated by $\overline{\NE}(S)_{A\geq 0}$ and 
$\sum_j \mathbb R _{\geq 0} [C_j]$.
This fact and (2) implies (1) (See the proof of \cite [Theorem 1.24] {KM}). 
(3) is also obvious. 
By changing $A$ to $A+\epsilon H\sim _{\bQ}K_X+B$, 
where $\epsilon$ is a small rational number, 
we can reduce it to the well-known klt case. 
\endproof

\ifx\undefined\bysame
\newcommand{\bysame}{\leavevmode\hbox to3em{\hrulefill}\,}
\fi

\end{document}